\documentclass[a4paper,12pt]{article}

\usepackage{amsmath}
\usepackage{amssymb}
\usepackage{theorem}

\begin{document}
\begin{center}
 \Large

 Appendix

 to ``Inconsistency of Inaccessibility''

 \medskip

 \large

  \textit{Alexander Kiselev}

 \smallskip

  aakiselev@yahoo.com
\end{center}

\renewcommand{\abstractname}{}
\begin{abstract}
{\footnotesize This paper is the concise addition to the foregoing
work ``Inconsistency of Inaccessibility'', containing the
presentation of main theorem proof (in\emph{ ZF}) about
inaccessible cardinals nonexistence. Here some refinement of this
presentation is set forth. Much attention is devoted to the
explicit and substantial development and cultivation of basic
ideas, serving as grounds for all main constructions and
reasonings.}

\end{abstract}

The proof of main theorem in its brief form is exposed in
specified work~\cite{Kiselev1}, one can also see it at arXiv site
~\cite{Kiselev1}. The present paper is its short complement using
previous notions and designations from~\cite{Kiselev1}. The proof
of main theorem has been published in 1997, its systematic
exposition has been published in 2000, and the most complete and
detailed form the proof of inaccessible cardinals nonexistence
have received in works ``Inaccessibility and Subinaccessibility",
Part I~\cite{Kiselev2} and Part II~\cite{Kiselev3} in 2008, 2010;
these two works one can see also at arXiv sites
~\cite{Kiselev2},~\cite{Kiselev3}.

Remind, that further matices and matrix functions reduced to \
$\chi^{\ast} $ \ are considered, the lower index \ $\chi^{\ast} $
\ will be omitted in notations, for instance \ $\alpha
_{\chi^{\ast}}^{\downarrow }$,\ \ $\alpha
_{\chi^{\ast}}^{\Downarrow }$,\ \ $S_{\chi^{\ast} f}$, \ \
$S_{\chi ^{\ast }\tau }$ will be denoted through \
$\alpha^{\downarrow }$,\ \ $\alpha^{\Downarrow }$,\ \ $S_{f}$, \ \
$S_{\tau }$  and so on.

Here it remains only to eliminate the case when \ $\delta =\chi
^{\ast +}$, \ that is when the function \ $S_{f}^{<\alpha^0}$ \ is
\ $\underline{\lessdot }$-nondecreasing up to \  $\chi ^{\ast +}$.
It means that the unrelativized function \ $S_f$, being \
$\underline{\vartriangleleft }$~-bounded by the fixed ordinal \
$Od \left ( S_{\tau^{\ast}} \right ) < \chi^{\ast +} $, \ loses
this boundedness after its relativization to \ $ \alpha^{0}$, that
is  to its transformation into \ $S_{f}^{<\alpha^0}$. \ One can
see that it happens because of deformation of the universe after
its relativization to \ $\alpha^0 $: \ because of losing
properties of subinaccessibility of all levels \ $\geq n$ \ by all
cardinals \ $\leq \gamma _{\tau _{0}}$ \  after their
relativization to \ $\alpha^0 $\ (Kiselev~\cite{Kiselev4}).
\\
So, many important properties of lower levels of the universe do
not extend up to relativizing cardinals, namely to prejump
cardinals of reduced matrices that are values of matrix functions.
And all the rest part of this work consists in correction of this
lack. In order to prevent this phenomenon we introduce special
cardinals named disseminators that are extending such properties
without distortion:
\\
\textbf{Definition 1 }\quad \ \emph{Let
    \quad \ \ $0<\alpha <\alpha _{1} \leq k,\quad X\subseteq \alpha, \quad
    X\neq \varnothing $.
\\
The ordinal \  $\alpha $ \  is named disseminator of level \ $m$
\  with data base \  $X$ \ below \  $\alpha _{1}$ \ iff for every
\ $\mathfrak{n} \in \omega_0$ \ and train \ $\overrightarrow{a}$ \
of  ordinals \  $\in X$ there holds \quad \ \ $U_{m}^{ \Sigma
\vartriangleleft \alpha _{1}} ( \mathfrak{n},
    \overrightarrow{a} ) \longrightarrow  U_{m}^{\Sigma
    \vartriangleleft \alpha } ( \mathfrak{n}, \overrightarrow{a} )$.
\\
The formula defining in \ $L_{k}$ \ the set of all such
disseminators is denoted by \ $ SIN_{m}^{<\alpha _{1}}\left[
X\right] ( \alpha  ) $; \ this set itself is also denoted by \
$SIN_{m}^{<\alpha _{1}}\left[ X \right] $, \ while its
disseminators -- by the common symbol \  $\delta ^{X}$ \  or,
briefly, by \ $\delta $, \ pointing to \ $\alpha_1$, \ \ $X$ in
the context. Everywhere further \emph{\emph{m=n+1}} will be
considered.}   \hspace*{\fill} $\dashv$
\\
Remind, the symbol \ $U_{n+1}^{\Sigma}$ \  here denotes the \
$\Sigma _{n+1}$-formula which is universal for the class of \
$\Sigma _{n+1}$-formulas but without any occurrences of the
constant \ $\underline{l}$. \  The  term ``disseminator'' is
justified by the fact that such an ordinal really extends \
$\Pi_{n+1}$-properties (containing constants from its base) from
lower levels of the universe up to relativizing cardinal \ $\alpha
_{1}$.
\\
The discussion of disseminators is more clear when\  data base \
$X$ \ is some ordinal \ $\rho$ \ (that is the set of all smaller
ordinals). In this case we shall write
    \ $SIN_{n+1}^{<\alpha _{1}}\left[ <\rho \right]$
    instead of \ $SIN_{n+1}^{<\alpha _{1}}
    \left[ X\right]$;
further \ $\rho$ is closed under the pair function.  Now one
should supply singular matrixes with disseminators:
\\
\textbf{Definition 2}  \emph{Let \
$\chi^{\ast}<\alpha<\alpha_1\leq k $ \ and \ $S$ \ be a matrix
reduced to cardinal \ $\chi^{\ast}$ \ and singular on a carrier \
$\alpha $.}
\\
\emph{1) \quad We name as disseminator for \ $S$ \  on \ $\alpha $
\ (or as disseminator for this carrier) of the level \  $n+1$ \
with data base \  $\rho $ \  every disseminator
 \ $\delta \in SIN_{n+1}^{<\alpha^{\Downarrow }}\left[ <\rho \right]
\cap SIN_{n}^{<\alpha^{\Downarrow }}$  below the prejump cardinal
\ $\alpha^{\Downarrow }$ \ after \ $\chi^{\ast}$.}
\\
\emph{2) \quad The matrix \ $S$ \ is called the disseminator
matrix or, more briefly, the \ $\delta $-matrix of the level \
$n+1$ \ admissible on the carrier \ $\alpha$ for \ $\gamma =
\gamma_{\tau}^{<\alpha_1}$ \ below \ $\alpha_1$ \ iff it possesses
some disseminator \ $ \delta < \gamma$ \ of the level \ $n+1$ \
with some base \ $\rho \le \chi^{\ast+}$ \ on this carrier such
that \ $S \vartriangleleft \rho$ \  (also below \ $\alpha_1$). }
\hspace*{\fill} $\dashv$

Now we are going to undertake the second approach  to the main
theorem proof idea. Remind that the simplest matrix function \
$S_{f}$ \
  possess the property of monotonicity, but the
direct proof of the required contradiction -- the proof of its
nonmonotonicity -- is hampered by the occasion: some essential
properties of lower levels of universe do not extend up to prejump
cardinals of matrices on their carriers that are values of matrix
functions.
\\
In order to destroy this obstacle we equipped such matrices with
corresponding disseminators and as the result the simplest matrix
functions are transformed to their more complicated forms, that is
to \ $\delta $-functions:
\\
\textbf{Definition 3} \quad \  \emph{Let \quad \ \ $\gamma <
\alpha <\alpha _{1}\leq k$. \quad 1) We denote through \
$\mathbf{K}_{n}^{\forall <\alpha _{1}}(\gamma ,\alpha )$ \ the
formula: }

\quad \quad \quad \quad \quad \quad \ \ $SIN_{n-1}^{<\alpha
_{1}}(\gamma )\wedge \forall \gamma ^{\prime }\leq \gamma \
(SIN_{n}^{<\alpha _{1}}(\gamma ^{\prime })\longrightarrow
SIN_{n}^{<\alpha }(\gamma ^{\prime }))~$.
\\
\emph{If this formula is fulfilled by constants \ $\gamma $, \
$\alpha $, \ $\alpha _{1}$, \ then we say that \ $\alpha $ \
conserves \ $SIN_{n}^{< \alpha _{1}}$-cardinals \ $\leq \gamma $ \
below \ $\alpha _{1}$.}
  \emph{If \ $S$ \ is a matrix on a carrier
\ $\alpha $ \ and its prejump cardinal \ $\alpha^{\Downarrow }$ \
after \ $\chi^{\ast}$ \ conserves these cardinals, then we also
say that \ $S$ \ on \ $\alpha $ \ conserves these cardinals below
\ $\alpha _{1}$.}
\\
 2)  \emph{We denote through \ $\mathbf{K}_{n+1}^{\exists }(\chi^{\ast} ,\delta
,\gamma ,\alpha ,\rho ,S)$ \ the \ $\Pi _{n-2}$-formula:}

\quad \quad \ $\sigma (\chi^{\ast} ,\alpha ,S)\wedge \chi^{\ast}
<\delta <\gamma <\alpha \wedge S\vartriangleleft \rho \leq
\chi^{\ast+}\wedge SIN_{n+1}^{<\alpha^{\Downarrow }}\left[ <\rho
\right] (\delta )\wedge SIN_{n}^{<\alpha^{\Downarrow }}(\delta )$.

\noindent \emph{Here, remind, the formula \ $\sigma(\chi^{\ast},
\alpha, S)$ \ means that \ $S$ \ is singular matrix on its carrier
\ $\alpha$ \ reduced to the cardinal \ $\chi^{\ast}$ ;\ \ $\delta$
\ is disseminator for \ $S$ \ on \ $\alpha$ \ with the base \
$\rho$ \ of the level \ $n+1$. \
 We denote through \ $\mathbf{K}^{<\alpha _{1}}(\chi^{\ast} ,\delta
,\gamma ,\alpha ,\rho ,S)$ \ the formula:}

\quad \quad \quad \quad \ $\mathbf{K}_{n}^{\forall <\alpha
_{1}}(\gamma ,\alpha^{\Downarrow })\wedge
\mathbf{K}_{n+1}^{\exists \vartriangleleft \alpha
_{1}}(\chi^{\ast} ,\delta ,\gamma ,\alpha ,\rho ,S)\wedge \alpha
<\alpha _{1}~$.
\\
3) \emph{If this formula is fulfilled by constants \
$\chi^{\ast}$, $\delta $, $\gamma $, $\alpha $, $\rho$, $S$,
$\alpha _{1}$, \ then we say that \ $\chi^{\ast} $, $\delta $,
$\alpha $, $\rho $, $S$ \ are strongly admissible for \ $\gamma $
\ below \ $\alpha _{1}$.
 If some of them are fixed or meant by the context, then
we say that others are strongly admissible for them (and for \
$\gamma $) \ below \ $\alpha _{1}$.\emph{} }
\\
\emph{4) \quad The matrix \ $S$ \ is called strongly disseminator
matrix or, briefly, \ $\delta $-matrix strongly admissible on the
carrier \ $\alpha $ \ for \ \mbox{$\gamma =\gamma _{\tau
}^{<\alpha _{1}}$} \ below \ $\alpha _{1}$, \ iff it possesses
some disseminator \ $\delta <\gamma $ \ with base \ $\rho $ \
strongly admissible for them (also below \ $\alpha _{1}$). In
every case of this kind \ $\delta $-matrix is denoted by the
common symbol \ $\delta S$ \ or \ $S$.}  \hspace*{\fill} $\dashv$
\\
\textbf{Definition 4 }\quad \ \emph{Let \ $\chi ^{\ast }<\alpha
_{1}$.
 \quad We call as the matrix \ $\delta $-function of the level \
$n$ \ below \ $\alpha _{1}$ \ reduced to \ $\chi ^{\ast }$ \ the \
function
 \quad \ \ $\delta S_{f}^{<\alpha _{1}}=(\delta S_{\tau }^{<\alpha
_{1}})_{\tau }~$
 \quad \  taking values for \ $\tau<k$:}

\quad \quad \quad \quad \ \ $\delta S_{\tau }^{<\alpha
_{1}}=\min_{\underline{\lessdot }} \bigl \{ S \vartriangleleft
\chi^{\ast +} : \exists \delta ,\alpha ,\rho <
\gamma_{\tau+1}^{<\alpha_1} ~ \mathbf{K}^{<\alpha _{1}}(\chi
^{\ast }, \delta ,\gamma _{\tau }^{<\alpha _{1}},\alpha ,\rho ,S)
\bigr \}$. \hspace*{\fill} $\dashv$
\\
Therefore such function has values that are singular strongly
admissible on their carriers \ $\delta $-matrices.

 The unrelativized function \ $\delta S_{f} $ \ really does exist
on the final subinterval of the inaccessible cardinal \ $k$ \ as
it shows
\\
\textbf{Lemma 5} \emph{(About \ $\delta $-function definiteness)
 There exists an ordinal \ $\delta <k$ \
such that \ $\delta S_{f}$ \ is defined on the set
 \quad \quad \ \ $T = \{\tau :\delta <\gamma _{\tau } < k \}$.
\\
The minimal of such ordinals \ $\delta $ \ is denoted by \ $\delta
^{\ast }$ \ and its index by
    \ $\tau_1^{\ast }$, so that
    \ $\delta^{\ast}=\gamma_{\tau_1^{\ast}}$}
 \hspace*{\fill} $\dashv$

 The monotonicity of this function is treated in the previous way: the function \ $\delta
S_{\tau}^{<\alpha_1}$ \  is called (totally) monotone iff \quad
\quad \ \ $\forall \tau^{\prime},
\tau^{\prime\prime}(\tau_1^{\ast} < \tau^{\prime}
    < \tau^{\prime\prime} < k \longrightarrow \delta
    S_{\tau^{\prime}} \; \underline{\lessdot } \; \delta S_{\tau^{\prime\prime}} )$.
\\

Let's discuss the situation which arises.  For some convenience
the suitable notation of ordinal intervals will be used for \
$\alpha _{1}<\alpha _{2}$:
 $\left[ \alpha _{1},\alpha _{2}\right[ =\alpha _{2}-\alpha _{1}$;
\quad $\left] \alpha _{1},\alpha _{2}\right[ =\alpha _{2}- (
\alpha _{1}+1 ) $ and so on (here ordinals \ $\alpha _{1}, \alpha
_{2}$ \ are sets of smaller ordinals).
\\
It was noted above, that the simplest matrix function \ $S_{f}$ \
is \ $ \underline{ \lessdot }$-monotone, but for every \
$\tau>\tau ^{\ast }$ \ the prejump cardinal \ $\alpha ^{\Downarrow
}$ \ of \ $S_{\tau} $ \ on its corresponding carrier \ $\alpha \in
\; ]\gamma_{\tau}, k [ $ \  \emph{do not conserve the
subinaccessibility of levels \ $\geq n$\ } of cardinals \ $\leq
\gamma _{\tau }$, \ and some other important properties of the
lower levels of the universe are destroyed when relativizing to \
$ \alpha ^{\Downarrow }$. \

In order to overcome these obstruction we have supplied the values
of this function, matrices \ $ \delta S_{\tau }$, \  by
disseminators of the level \ $n+1$ \ and required the conservation
of the subinaccessibility of the level \ $n$ \ for cardinals \
$\leq \gamma _{\tau }$, \ that is we passed to the \ $\delta
$-function \ $\delta S_{f}$. And these disseminators by
definitions 2, 3 really do overcome these obstructions and now
basic $\Pi_{n+1}$-properties and all $\Pi_{n}$-properties of the
lower levels of the universe are extending up to prejump cardinals
of matrixes \ $ \delta S_{\tau }$ \ carriers after \ $\chi ^{\ast
}$. \

 But now it involves the new complication: now one can see, that after
this modification such \ $\delta $-function loses its  property of
(total) monotonicity on $[\tau_1^{\ast}, k[\;$, \ and just due to
the fact that prejump cardinals \ $\alpha^{\Downarrow } $ \ of \
$\delta$-matrices carriers \ $\alpha$, \ vice versa,
\textit{generate the subinaccessibility of the level} \ $n$ \ of
some cardinals \ $\leq \gamma _{\tau }$ \ that become
subinaccessible (relatively to \ $\alpha^{\Downarrow} $), \ not
being those in the universe (Kiselev~\cite{Kiselev5}).

So, here we come to the third and final approach to the main
theorem proof idea:
\\
 \textit{The following requirements should
be imposed on \ $\delta$-matrices \ $ \delta S_{\tau }$ \
successively in the course of matrix function defining:
\\
1) they must possess the property of ``autoexorcizivity'', that is
of self-exclusion in situations when \ $ S= \delta S_{\tau }$ \ is
the first value of \ $\delta$-function monotonicity violation in
its preceding part; the matrices with this property (of {\sl
``unit characteristic''}) will have the priority over other
matrices (of {\sl ``zero characteristic''} respectively) in the
course of defining the matrix function;
\newline 2) more requirements should be imposed on matrices of {\sl
zero} characteristic, preventing their forming: their disseminator
data bases must increase substantially, and must be even disabled
(``suppressed''), when the preceding part of matrix function,
which has already been defined, contains monotonicity violation,
in order to correct this fault -- the using of matrices of {\sl
zero} characteristic;
\\
on this grounds \ $\delta$-matrix functions should receive
inconsistent properties of monotonicity and nonmonotonicity
simultaneously. }

 Obviously, all these considerations require the
recursive definition of matrix functions, setting there values
depending on specified properties of their preceding values.
\\
This resulting recursive definition follows definition 4, but in
recursive form and the formula \ $\mathbf{K}^{<\alpha _{1}}$ \
transforms in its corresponding version which we denote through \
$\alpha \mathbf{K}^{\ast <\alpha _{1}}$. \ As the result there
arises matrix functions which we denote through \ $\alpha
S_{f}^{<\alpha _{1}}$ and name \ $\alpha $-functions below \ $
\alpha _{1}$.
\\
 One can see that these functions
  possess many basic properties of simplest matrix functions and \ $\delta
$-functions, for instance, analogous to their definiteness
properties for unrelativized  \ $\alpha $-function \ $ \alpha
S_{f}$: there exists the cardinal \ $\delta <k$ \ such that
 \quad \ \ $\{ \tau: \delta < \gamma _{\tau } < k \} \subseteq  dom ( \alpha
S_{f} )$;
 \quad \ more precisely:
\\
\textbf{Lemma 6} \emph{(About \ $\alpha $-function definitness)
 \quad \  There exist cardinals \ $\delta < \gamma <
k$ \ such that for every \ \mbox{$SIN_n$-}car\-di\-nal \ $\alpha_1
> \gamma$, $\alpha_1 < k$ \ limit for \ $SIN_n \cap \alpha_1$ \ the function \ $\alpha
S_f^{<\alpha_1}$ \ is defined on the nonemty set
    \quad \ \ $T^{\alpha_1} = \{\tau: \delta < \gamma_{\tau}^{<\alpha_1} < \alpha_1\}$.}
\\
\emph{The minimal of such cardinals \ $\delta$ \ is denoted by \
$\alpha \delta ^{\ast }$ \  and its index by}
    \quad \ \ $\alpha \tau_1^{\ast} = \tau(\alpha \delta^{\ast})$,
  \quad \ \emph{so that}
    \quad \ \ $ \alpha \delta^{\ast} = \gamma_{\alpha
    \tau_1^{\ast}}$.\      \hspace*{\fill} $\dashv$

\hspace*{1em} Now the monotonicity of \ $\alpha$-function is
considered; here this notion is used again in the previous sense:
 the function \ $\alpha
S_{f}^{<\alpha_{1}}$ \ is called monotone on interval \ $\left[
\tau _{1},\tau _{2}\right[ $ \  below \ $\alpha_1$ \ iff  \
$\left] \tau _{1},\tau _{2} \right[ \subseteq dom ( \alpha
S_{f}^{<\alpha _{1}} ) $ \ and \quad \  $\forall \tau ^{\prime
},\tau ^{\prime \prime }  ( \tau _{1} < \tau^{\prime }<\tau
^{\prime \prime } < \tau _{2}\longrightarrow \alpha S_{\tau
^{\prime }}^{<\alpha _{1}} \underline{\lessdot }\alpha S_{\tau
^{\prime \prime }}^{<\alpha _{1}} )~$. As we shall see, this
property is rather strong; in particular, any interval \ $[
\tau_1, \tau_2 [$ \ of its monotonicity can not be ``too long'' --
that is the corresponding interval \ $] \gamma_{\tau_1},
\gamma_{\tau_2} [$ \ can not contain any \
\mbox{$SIN_n$-cardinals}; therefore \ $\alpha S_{f}^{<\alpha_{1}}$
\ doesn't possess this property in total sense:
\\
 \textbf{Theorem 1.}
 \emph{ Let \medskip }
\emph{ \ $\alpha S_f^{<\alpha_1}$ \ be monotone on \ $[\tau_1,
\tau_2[ $ \ below \ $\alpha_1$ \ and } \emph{ \ $\tau_1 = \min \{
\tau: \; ]\tau, \tau_2[ \; \subseteq$  \ $ \subseteq dom (\alpha
S_f^{<\alpha_1}) \} $,}  \emph{ then} \quad \quad \quad \quad \ \
$]\gamma_{\tau_1}^{< \alpha_1}, \gamma_{\tau_2}^{< \alpha_1} [ ~
\cap ~ SIN_n^{<\alpha_1} = \varnothing$. \hspace*{\fill} $\dashv$

The contradiction, which proves the main theorem, is the
following:
\\
On one hand, by lemma 6  function \ $\alpha S_f^{<\alpha_1}$ \ is
defined on nonempty set
 \ $T^{\alpha_1} = \{\tau:\alpha\delta^{\ast}<$ \ $ <\gamma_{\tau}\ < \alpha_1 \}$
for every sufficiently great cardinal \ $\alpha_1 \in SIN_n$. Its
monotonicity on this set is excluded by theorem 1.
\\
But on the other hand, this monotonicity is ensured by the
following theorem for every \ $SIN_n$-cardinal \
$\alpha_1>\alpha\delta^{\ast}$ \ of sufficiently great cofinality:
\\
\textbf{Theorem 2.} \emph{Let function \ $\alpha S_f^{<\alpha_1}$\
be defined on nonempty set} \ $T^{\alpha_1}
=\{\tau:\gamma_{\tau_1}^{<\alpha_1} < \gamma_{\tau}^{<\alpha_1} <
\alpha_1 \}$
\\
\emph{for }\ $\alpha_1 < k $ \ such that
 \quad \  (i) \ $\tau_1=\min\{\tau:\forall\tau^{\prime}
(\gamma_{\tau}^{< \alpha_1} <
\gamma_{\tau^{\prime}}^{<\alpha_1}\longrightarrow\tau^{\prime}\in
dom(\alpha S_f^{<\alpha_1}))\},$ \
\\
 (ii) \ $\sup SIN_n^{<\alpha_1}=\alpha_1;$
 \quad \ (iii) \ $cf(\alpha_1)\geq\chi^{\ast +}.$
\\
\noindent \emph{Then \ $\alpha S_f^{<\alpha_1}$ \ is monotone on
this set below \ $\alpha_1$:}
    \quad \ \ $ \forall \tau_1, \tau_2 \in T^{\alpha_1} \bigl( \tau_1 <
    \tau_2 \rightarrow \alpha S_{\tau_1}^{<\alpha_1}
    \underline{\lessdot} \alpha S_{\tau_2}^{<\alpha_1} \bigr)$.
\\
\noindent \textit{Proof}  is the following in general. The
reasoning will be carried out by the induction on the cardinal \
$\alpha_1$.
 Suppose, that this theorem fails and the cardinal \
$\alpha_1^\ast$ \ is \textit{minimal} breaking this theorem, that
is the function \ $\alpha S^{<\alpha_1^\ast}_f$ \ is nonmonotone
on  the nonempty set
 \ $T^{\alpha_1^\ast} = \bigl \{ \tau:$
\ $ :\gamma_{\tau_1^\ast}^{<\alpha_1^\ast} <
\gamma_{\tau}^{<\alpha_1^\ast} < \alpha_1^{\ast} \bigr \}$ with
specified properties $(i)$--$(iii)$ for some \ $\tau_1^{\ast}$, \
so that the \textit{inductive hypothesis} is accepted:
\\
for every \ $\alpha_1 < \alpha_1^{\ast}$ \ the function \ $\alpha
S_f^{<\alpha_1}$ \ is monotone on the nonempty set \
$T^{\alpha_1}$ \ with properties $(i)$--$(iii)$.
\\
It follows straight from theorem 1, that this \ $\alpha_1^{\ast}$
\ is simply the \textit{minimal} cardinal \ $\alpha_1$, \
providing the existence of such nonempty set \ $T^{\alpha_1}$ \
because for every such \ $\alpha_1 < \alpha_1^{\ast}$ \ the
function \ $\alpha S_f^{<\alpha_1}$ \ on \ $T^{\alpha_1}$ \ is
\textit{nonmonotone} by theorem~1 and at the same time is
\textit{monotone} by the minimality of \ $\alpha_1^{\ast}$.

 But in this situation all possible violations
of function \ $\alpha S^{<\alpha_1^\ast}_f$ \  monotonicity on
specified set
 \ $T^{\alpha_1^\ast}$ draws the existence of some nonempty set
\ $T^{\alpha_1}=\{\tau:\gamma_{\tau_1}^{<\alpha_1} <
\gamma_{\tau}^{<\alpha_1} < \alpha_1 \}$ with properties
$(i)$--$(iii)$ for lesser \ $\alpha_1 < \alpha_1^{\ast}$. \ So,
the function \ $\alpha S^{<\alpha_1^\ast}_f$,\ nevertheless,  is
monotone on  the set
 \ $T^{\alpha_1^\ast}$.              \hspace*{\fill} $\dashv$

 Now let us sum up.
All the reasoning was conducted in the system

\quad \quad \quad \quad \quad \quad \ \ $ ZF + \exists k ~ (k
\mbox{ is weakly inaccessible cardinal})$,\
\\
where it  was considered the countable standard model
    \quad \ \ $ \mathfrak{M} = (L_{\chi^0}, \in, =)$
\quad \ of some finite part of the theory
  \quad \ \ $ ZF + V = L + \exists k ~ (k \mbox{ is weakly inaccessible
  cardinal})$.
\\
In the model \quad \ \ $(L_{k}, \in, =)$ \quad \ matrix \ $\alpha
$-functions were considered; such function \ $\alpha
S_f^{<\alpha_1}$ \ is defined on any nonempty set \
$T^{\alpha_1}$, \ which exist for any sufficiently great cardinal
\ $\alpha_1 < k$, $\alpha_1 \in SIN_n$. \
\\
 It provides the final
contradiction: take any \ $SIN_n$-cardinal \ $\alpha_1 > \alpha
\delta^{\ast}$ \ limit for \ $SIN_n \cap \alpha_1$ \ of the
cofinality \ $cf(\alpha_1) \ge \chi^{\ast +}$ \ providing such
nonempty set \ $T^{\alpha_1}$ \ with properties $(i)$--$(iii)$
from theorem~2; such set exist due to lemma 6; then the function \
$\alpha S_f^{<\alpha_1}$ \ is nonmonotone on this \ $T^{\alpha_1}$
\ by theorem~1 and at the same time is monotone on this set by
theorem~2. This contradiction ends the main theorem proof.
\\

\end{document}